\documentclass[12pt]{article}
\usepackage{amsfonts}
\usepackage{graphicx}
\begin{document}

\newtheorem{theorem}{Theorem}[section]
\newtheorem{prop}[theorem]{Proposition}
\newtheorem{cor}[theorem]{Corollary}
\newtheorem{lem}[theorem]{Lemma}
\newtheorem{definition}[theorem]{Definition}
\newtheorem{ex}[theorem]{Example}
\newtheorem{no}[theorem]{Note}
\newtheorem{unnumber}{}
\renewcommand{\theunnumber}{\relax}
\newtheorem{prepf}[unnumber]{Proof}
\newenvironment{pf}{\prepf\rm}{\endprepf}
\newcommand{\qed}{\qquad$\Box$}

\title{Catalan tree \\ \& \\
 Parity of some Sequences which are related to Catalan numbers }
\author{Volkan Yildiz\\\\
Department of Mathematics\\
King's College London\\
Strand, London WC2R 2LS\\
\texttt{volkan.yildiz@kcl.ac.uk}}
\date{May 24, 2011}
\maketitle

\begin{abstract} 
In this paper we determine the parity of some sequences which are related to Catalan numbers. 
Also we  introduce a combinatorical object called, ``Catalan tree'', and discuss its properties.
\end{abstract}

{\footnotesize

}

\section{Introduction}

Throughout this paper, for brevity, we represent the set of even counting numbers by 
the capital letter $E$, the set of odd counting numbers by the capital letter $O$, 
and the set of natural numbers, $\{1,2,3,4,...\}$, by $\mathbb{N}$.\\

In this paper we first define the term {\em Catalan tree} and then study its combinatorical properties. 
Later we explore the parity of some sequences which are related to Catalan numbers. 

The Catalan numbers are an infinite sequence of integers 1, 1, 2, 5, 14, 42, 132, 429, 1430, 4862, $\ldots$.
They are defined by the following recurrence relation:
\begin{equation}\label{e:1}
 C_n=\sum_{i=1}^{n-1}C_i C_{n-i}, \textit{ with } C_0=0, \; C_1=1.
\end{equation}
It also has the following explicit formula $C_n=\frac{1}{n}{{2n-2}\choose{n-1}}$ 
and for large $n$,  $C_n$ behaves like $\frac{2^{2n}}{\sqrt{\pi n^3}}$, see~\cite{P1}.
\\\\
The Catalan number appears in many areas of mathematics  
just like the Fibonacci number. Specifically they are related in 
combinatorical settings such as trees, lattice paths, partitions, (see~\cite{C})
and even within propositional logic, (see~\cite{P1}). 
Here we go one step further and define what a Catalan tree is.
\begin{definition}
The nth {\bf Catalan tree}, $A_n$, is a combinatorical object, 
characterized by one root, $(n-1)$  main-branches, and $C_n$  sub-branches.  
Where each main-branch gives rise to a number of sub-branches, and 
the number of these sub-branches is determined by the additive 
partition of the corresponding Catalan number,
 as determined by the recurrence relation (\ref{e:1}).
\end{definition}
The tree $A_n$ can be represented symbolically as follows:
\begin{eqnarray*}
C_n\textit{ sub-branches: } & (C_1C_{n-1},C_2C_{n-2},\ldots,C_{n-2}C_2,C_{n-1}C_1) \\
(n-1) \textit{ main-branches: } &  (1,1,\ldots,1,1) \\
\textit{one root: } & (1) 
\end{eqnarray*}
Note that the main-branches and the sub-branches exist iff $n>1$. 
Here is an example: The Catalan tree $A_4$, has one root, (1), followed by 3 main-branches,
 $(1,1,1)$, and each main-branch gives rise to $(2,1,2)$ sub-branches respectively. 
Also this combinatorical object can be represented by a graph:

\includegraphics[scale=0.29]{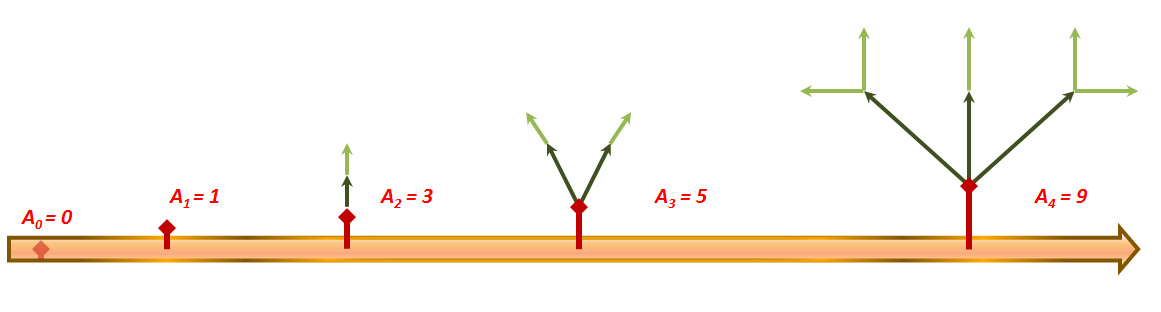}

The diagram above shows the first five stages of the Catalan tree, $A_n$, 
where $a_n$ is defined in Proposition~\ref{p:1}.

\begin{prop}\label{p:1} For $n>1$, let $A_n$ denote the nth Catalan tree, and 
let $a_n$ denote the number of components of $A_n$. Then 
\[
a_n=C_n+n,\;\; \;\;\textit{ with }\;\; a_0=0, a_1=1 .  
\]
\end{prop}

\begin{pf}
By definition there are 1 root, $(n-1)$ main branches and $C_n$ sub-branches. 
Therefore $a_n = C_n + n$, for $n>1$, and $a_0=0$, $a_1=1$.\qed
\end{pf}

Using Proposition~\ref{p:1}, it is straightforward to calculate the values of $a_n$. 
The table below illustrates this up to $n=10$.

\[
\begin{array}{|l|c|c|c|c|c|c|c|c|c|c|c|}
\hline n &0 & 1& 2 & 3 & 4 & 5 & 6 & 7 & 8 & 9 & 10\\
\hline a_n & 0 & 1 & 3 & 5 & 9 & 19 & 48 & 139 & 437 & 1439 & 4872\\
\hline
\end{array}
\]

Let $A(x)$, $C(x)$ and $N(x)$ denote the generating functions for $a_n$, $C_n$ and $n$, respectively. 
Thus 
\[
A(x)=\sum_{n\geq 1} a_n x^n, \; C(x)= \sum_{n\geq 1} C_n x^n = \frac{1}{2}(1-\sqrt{1-4x}), \; N(x)=\sum_{n\geq 1} nx^n = \frac{x}{(1-x)^2} .
\] 
\begin{prop} The generating function for the sequence $\{a_n\}_{n\geq 1}$ is given by
\[
A(x)= \frac{2x^2(2-x)+(1-x)^2(1-\sqrt{1-4x})}{2(1-x)^2}.
\]
\end{prop}

\begin{cor} For $n>1$, the explicit formula for $a_n$ is given by,
\[
a_n = \frac{1}{n}{{2n-2}\choose{n-1}} + n. 
\]
\end{cor}
The following result is a consequence of the asymptotic behavior of Catalan numbers.  
\begin{cor}
For large $n$, we have the asymptotic formula
\[
a_n \sim \frac{2^{2n}+n^2 \sqrt{\pi n}}{\sqrt{\pi n^3}}.
\]
\end{cor}
\section{Parity of related sequences}
In this section we determine the parity of the sequences that we have discussed in \cite{P1}, 
and as well as the parity of the sequences which are related to the sequences in \cite{P1}.

\begin{no}
The following Lemma \ref{l:1} has  been proven by number of other authors, (see~\cite{O}, (1986)),  
(see~\cite[page 330]{C}, (2004)) and (see~\cite{W}, (2008)). But, they took $n$ to be a Mersenne number, 
this is due to the fact that the Catalan numbers were shifted by one term forward in their work.
\end{no}
From the Segner's recurrence relation, $C_n$ can be expressed as a piecewise function, 
with respect to the parity of $n$, (see~\cite[page, 329]{C}).

\[C_n=\cases{2(C_1C_{n-1}+C_2C_{n-2}+\ldots+C_{\frac{n-1}{2}}C_{\frac{n+1}{2}}) \;\;& if $n\in O$,\cr\cr
2(C_1C_{n-1}+C_2C_{n-2}+\ldots+C_{\frac{n-2}{2}}C_{\frac{n+2}{2}})+C_{\frac{n}{2}}^2 \;\; & if $n\in E$.\cr}\]

\begin{lem}[Parity of $C_n$] \label{l:1}
\[ C_n\in O \Longleftrightarrow n=2^i, \textit{ where } i\in\mathbb{N} . \]
\end{lem}
\begin{pf}
For $n\geq 2$ 
\[
C_n\in O \Longleftrightarrow C_{\frac{n}{2}}^2 \in O  \Longleftrightarrow C_{\frac{n}{2}}\in O \Longleftrightarrow n=2^i \;\;\forall i\in \mathbb{N}. 
\]
Note that $C_1=1\in O$. \qed
\end{pf}

\begin{cor}[Parity of $a_n$]
\[
a_n \in O \Longleftrightarrow  n=2^i \textit{ or } n\in O,
\] (and $a_n\in E$ iff $n\in E$ and $n\not= 2^i$), for $i\in \mathbb{N}$.
\end{cor}
\begin{pf} Let the symbols $\wedge$, and $\vee$ denote the connectives `and' and `or' respectively. Then

\begin{eqnarray*}
a_n = (C_n+n) \in O & \Longleftrightarrow & (C_n \in O \wedge n \in E) \vee (C_n\in E \wedge n\in O)\\
& \Longleftrightarrow & (n=2^i \wedge n \in E) \vee (n\not= 2^i \wedge n\in O)\\
& \Longleftrightarrow & (n=2^i) \vee (n\in O) .
\end{eqnarray*}
\qed
\end{pf}

\begin{theorem} \label{t:f}
Let $f_n$ be the number of rows with the value ``false'' in the truth tables 
of all bracketed formulea with $n$ distinct propositions $p_1,\ldots,p_n$ 
connected by the binary connective of implication. Then in~\cite{P1} 
we have shown that the following results are true:
\begin{equation}\label{e:2}
f_n =\sum_{i=1}^{n-1} (2^iC_i-f_i)f_{n-i}, \; \textit{ with } \; f_1=1
\end{equation}
and for large $n$, $\;f_n \sim \Bigg(\frac{3-\sqrt{3}}{6}\Bigg)\frac{2^{3n-2}}{\sqrt{\pi n^3}}$.
\end{theorem}
Using Theorem~\ref{t:f}, we get the following triangular table. 
Where the left hand side column represents the sum of the corresponding row.
\begin{center}
\begin{tabular}{rccccccccc}
$f_2$:&    &    &    &    &  1\\\noalign{\smallskip\smallskip}
$f_3$:&    &    &    &  1 &    & 3  \\\noalign{\smallskip\smallskip}
$f_4$:&    &    &  4 &    &  3 &    &  12 \\\noalign{\smallskip\smallskip}
$f_5$:&    &  19 &    &  12 &    &  12&    &  61 \\\noalign{\smallskip\smallskip}
$f_6$:&  104 &    &  57 &    &  48 &    &  61 &    &  344\\\noalign{\smallskip\smallskip}
\end{tabular}
\end{center}
\begin{theorem}[Parity of $f_n$]\label{f:c} The sequence $\{f_n\}_{n\geq 1}$ preserves the parity of $C_n$.
\end{theorem}

\begin{pf} If an additive partition of $f_n$, (which is determined by the recurrence 
relation~(\ref{e:2})), is odd, then it comes as a pair; i.e.  
\[
(2^iC_i-f_i)f_{n-i} \in O \Longleftrightarrow f_i, f_{n-i}\in O \Longleftrightarrow (2^{n-i}C_{n-i}-f_{n-i})f_i\in O .
\]
Hence, $\bigg((2^iC_i-f_i)f_{n-i}+(2^{n-i}C_{n-i}-f_{n-i})f_i \bigg) \in E.$

Thus, $f_n$ can be expressed as a piecewise function depending on the parity of $n$:
\[f_n=\cases{\sum_{i=1}^{\frac{n-1}{2}} ((2^iC_i-f_i)f_{n-i}+(2^{n-i}C_{n-i}-f_{n-i})f_i ) \;\;& if $n\in O$,\cr\cr
\bigg(\sum_{i=1}^{\frac{n-2}{2}} ((2^iC_i-f_i)f_{n-i}+(2^{n-i}C_{n-i}-f_{n-i})f_i )\bigg)+ (2^{\frac{n}{2}}C_{\frac{n}{2}}-f_{\frac{n}{2}})f_{\frac{n}{2}} \;\; & if $n\in E$.\cr}\]
Finally,  
\[
f_n\in O \Longleftrightarrow  (2^{\frac{n}{2}}C_{\frac{n}{2}}-f_{\frac{n}{2}})f_{\frac{n}{2}} \in O \Longleftrightarrow f_\frac{n}{2} \in O  \Longleftrightarrow  n=2^i, \;\; \forall i\in\mathbb{N}.
\]
Note that $f_1=1\in O$. \qed
\end{pf}

\begin{theorem} \label{t:t}
Let $g_n$ be the total number of rows in all truth tables for
bracketed implications with $n$ distinct variables. 
Let $t_n$ be the number of rows with the value ``true'' in the truth tables 
of all bracketed formulea with $n$ distinct propositions $p_1,\ldots,p_n$ 
connected by the binary connective of implication. Then in~\cite{P1} 
we have shown that the following results are true:
\[
t_n =g_n-f_n, \; \textit{ with } \; t_0=0
\]
and for large $n$, $\;t_n \sim \Bigg(\frac{3+\sqrt{3}}{6}\Bigg)\frac{2^{3n-2}}{\sqrt{\pi n^3}}$.
\end{theorem}

\begin{prop}[Parity of $t_n$]
The sequence $\{t_n\}_{n\geq 1}$ preserves the parity of $C_n$.
\end{prop}
\begin{pf}Since 
\[
t_n=g_n-f_n=2^n C_n-f_n, \textit{ with } n\geq 1 
\]
The sequence $\{g_n\}_{n\geq 1}$ is always even, and the sequence 
$\{f_n\}_{n\geq 1}$ preserves the  parity of $C_n$ by Theorem~\ref{f:c}. 
Therefore the sequence $\{t_n\}_{n\geq 1}$ preserves the parity of $C_n$.
\qed
\end{pf}

\section{A fruitful tree}

\begin{definition}
The Catalan tree $A_n$ is {\bf fruitful} iff each sub-branch of $A_n$ has fruits. 
We denote this new tree by $A_n(\mu_i)$, 
where $\{\mu_i\}_{i\geq 1}$ is the corresponding fruit sequence.
\end{definition}

\begin{ex}

Let $\{f_n\}_{n\geq 1}$ be the corresponding fruit sequence for the Catalan tree $A_n$. 
Then $A_n(f_n)$ has the following symbolic representation, 
{\small
\begin{eqnarray*}
&((2^1C_1-f_1)f_{n-1}, (2^2C_2-f_2)f_{n-2}, \ldots, (2^{n-2}C_{n-2}-f_{n-2})C_2, (2^{n-1}C_{n-1}-f_{n-1})f_1) \\
& (C_1C_{n-1},C_2C_{n-2},\ldots,C_{n-2}C_2,C_{n-1}C_1) \\
&  (1,1,\ldots,1,1) \\
& (1).
\end{eqnarray*}
}
\end{ex}

\begin{ex}
More concretely, $A_5(f_n)$ has the following symbolic representation:
\begin{eqnarray*}
& ((5,7,3,3,1), (5,7), (9,3), (15,13,13,9,11)) \\
& (5,2,2,5)  \\
& (1,1,1,1) \\
& (1) .
\end{eqnarray*}

The diagram below shows $A_5(f_n)$ as a graph:

\includegraphics[scale=0.30]{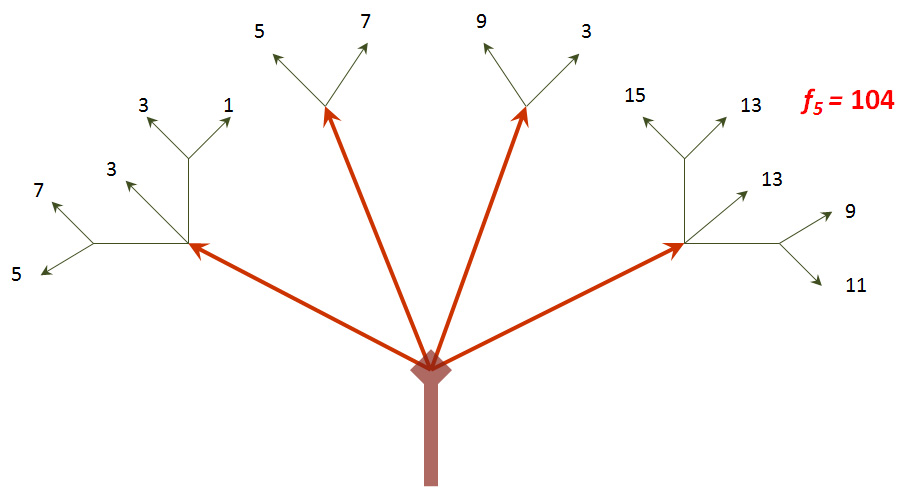}
\end{ex}

\begin{ex}
Let $\{t_n\}_{n\geq 1}$ be the corresponding fruit sequence for the Catalan tree $A_n$. 
Then $A_n(t_n)$  has the following symbolic representation, 
{\small
\begin{eqnarray*}
&(2^n-((2^1C_1-f_1)f_{n-1}), \ldots, (2^n-(2^{n-1}C_{n-1}-f_{n-1})f_1)) \\
& (C_1C_{n-1},C_2C_{n-2},\ldots,C_{n-2}C_2,C_{n-1}C_1) \\
&  (1,1,\ldots,1,1) \\
& (1).
\end{eqnarray*}
}
\end{ex}

\begin{prop}\label{p:s}  For $n>1$, let $a_n(f_n)$ and $a_n(t_n)$ be the total number of components of 
the fruitful trees $A_n(f_n)$ and $A_n(t_n)$ respectively. Then  
\[
a_n(f_n) = f_n + C_n + n, \; \textit{ and } \; a_n(t_n) = t_n + C_n + n .
\]
\end{prop}

Using Proposition~\ref{p:s}, it is straightforward to calculate the values of $a_n(f_n)$ , and $a_n(t_n)$. 
The table below illustrates this up to $n=10$.

\[
\begin{array}{|l|c|c|c|c|c|c|c|c|c|c|c|}
\hline n & 0 & 1& 2 & 3 & 4 & 5 & 6 & 7 & 8 & 9 & 10\\
\hline a_n(f_n)  & 0 & 2 & 4 & 9 & 28 & 123 & 662 & 3955 & 25032 & 164335 & 1106794 \\
\hline a_n(t_n)  & 0 & 2 & 6 & 17 & 70 & 363 & 2122 & 13219 & 85666 & 570703 & 3881638 \\
\hline
\end{array}
\]

\begin{cor} For $n\geq 1$, $a_n(f_n)$, and $a_n(t_n)$ are odd iff $n\in O$. \end{cor}
\begin{pf}
Since, 
\[
f_n, t_n \in O \Longleftrightarrow n=2^i, \textit{ and } a_n=(C_n+ n) \in O \Longleftrightarrow 
n=2^i \textit{ or } n\in O
\]
then $\; a_n(f_n), a_n(t_n)\in O \Longleftrightarrow n\in O$.	
 \qed
\end{pf}


\begin{thebibliography}{9}
\bibitem{P1}
P. J. Cameron and V. Yildiz,
\textit{Counting false entries in truth tables of bracketed formulae connected by implication},
Submitted to JIS, Preprint, 14-July-2010, ( arxiv.org/abs/1106.4443 ).

\bibitem{O}
\"{O}. E\~gecio\~glu,
\textit{ The parity of the Catalan numbers via lattice paths},
 Fibonacci Quart. 21 (1983) 65-66.

\bibitem{W}
 K.Q. Ji and H.S. Wilf, 
\textit{ Extreme Palindromes},
American Mathematical Monthly, 2008, VOL 115; NUMB 5, pages 447-450. 

\bibitem{C}
T. Koshy,
\textit{Catalan Numbers with Applications},
Oxford University Press, 2009.


\end{thebibliography}
\end{document}